\documentclass{amsart}
\usepackage{amsmath,amsfonts,amssymb,amsthm}
\usepackage{geometry}
\usepackage{hyperref}
\usepackage{graphicx}
\usepackage{tikz}
\usepackage{tikz-cd}
\usetikzlibrary{decorations.pathmorphing, arrows.meta}
\usepackage{verbatim}  

 \usepackage{amsfonts}
\usepackage{amsmath}
\usepackage{amssymb}
\usepackage{amstext}
\usepackage{amsthm}
\usepackage{epigraph}
\usepackage[percent]{overpic}

\usepackage{xypic}

\newtheorem{Theorem}{Theorem}[section]

\theoremstyle{definition}
\newtheorem{Definition}[Theorem]{Definition}
\newtheorem{Remark}[Theorem]{Remark}
\newtheorem{Example}[Theorem]{Example}
\usepackage{chngcntr}
\counterwithin{figure}{section}

\begin{document}

\title{A physical model of the Reeb foliation}
\author{G.~Bande}
\address{Dipartimento di Matematica e Informatica, Universit\`a degli Studi di Cagliari, Via Ospedale
72, 09124 Cagliari, Italy}
\email{gbande{\char'100}unica.it}
\author{G. Franzoni}
\address{Dipartimento di Matematica e Informatica, Universit\`a degli Studi di Cagliari, Via Ospedale
72, 09124 Cagliari, Italy}
\email{gregorio.franzoni{\char'100}unica.it}
\date{\today}
\begin{abstract}
In this paper, after explaining some basic aspects of the modern theory of foliations with the aim of describing the celebrated Reeb foliation, we propose the first construction of a comprehensive physical model of it. The construction of the model is obtained by an implementation of  geometric methods for 3D printing.
\end{abstract}
\maketitle

\section*{Introduction}
Foliations are an important structure in differential geometry and topology. They allow the decomposition of a differentiable manifold into a union of submanifolds, called \emph{leaves}, that fit together in a smooth and regular manner.
The concept of foliation can be traced back to the 19th century. However, the modern theory of foliations was established in the 20th century, primarily through the contributions of Georges Reeb and Charles Ehresmann. In particular, Georges Reeb was able to construct a foliation on the $3$-dimensional sphere such that all the leaves are diffeomorphic to $\mathbb{R}^2$ except one which is diffeomorphic to the $2$-dimensional torus. It was the first example of foliation having non-diffeomorphic leaves and is called after him the {\it Reeb foliation}.

Topologically, the $3$-sphere is obtained by gluing two solid tori along their boundary, which is a torus. The global foliation is then obtained by foliating the interior of each solid torus by surfaces diffeomorphic to planes, in such a way that these surfaces accumulate to the boundary torus, which also is a leaf of the solid torus. The solid tori endowed with these foliations are the so-called \emph{Reeb components}. The resulting foliation of the sphere has just one compact leaf (given by the glued boundary tori) and the other leaves are in fact diffeomorphic to $\mathbb{R}^2$.

We refer to the nice survey of André Haefliger \cite{H} for more explanations of the seminal work of Reeb and Ehresmann, the contributions of Haefliger himself and the fundamental result of Novikov \cite{N} which states, in particular,  that every foliation of the $3$-sphere must contain a Reeb component.

\begin{figure}[h]
\caption{Two rendered images of the Reeb foliation in $\mathbb{R}^3$, showing four leaves belonging to the ``external" Reeb component (the blu ones). The white surface is the compact leaf.}\label{first rendering 1}
\includegraphics[width=0.8\textwidth]{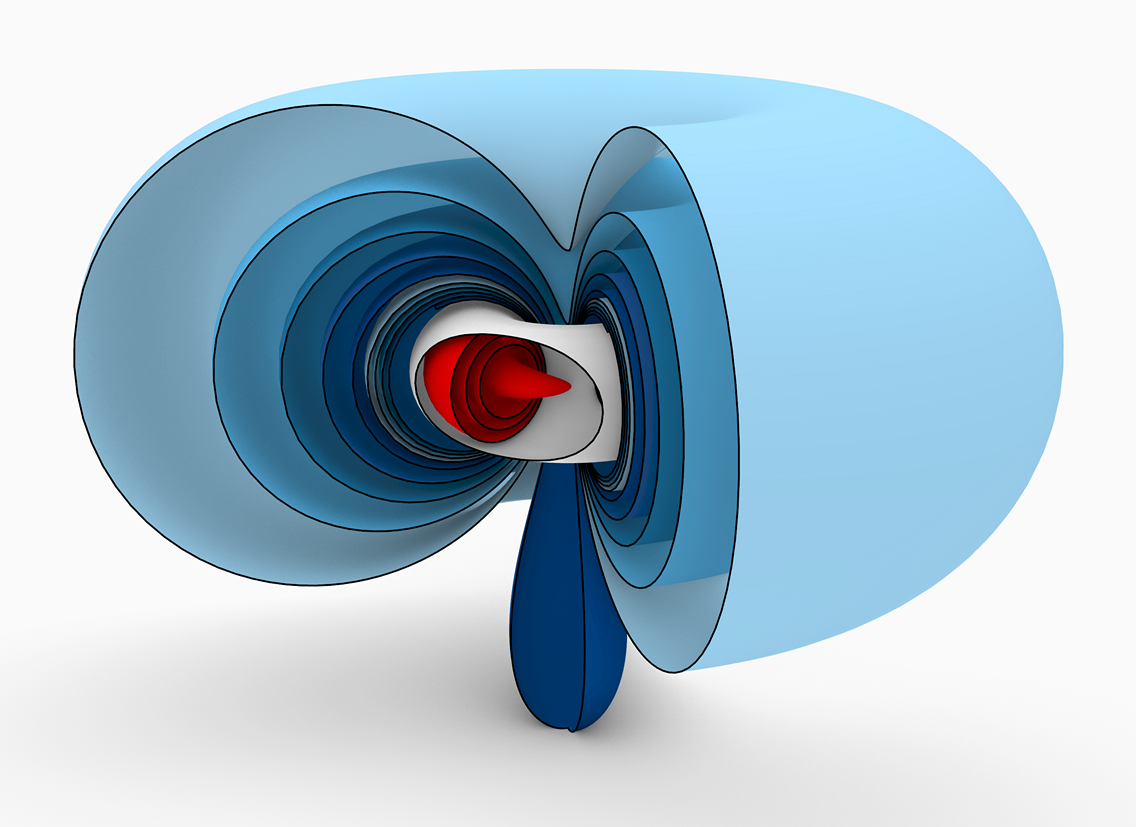}

\includegraphics[width=0.8\textwidth]{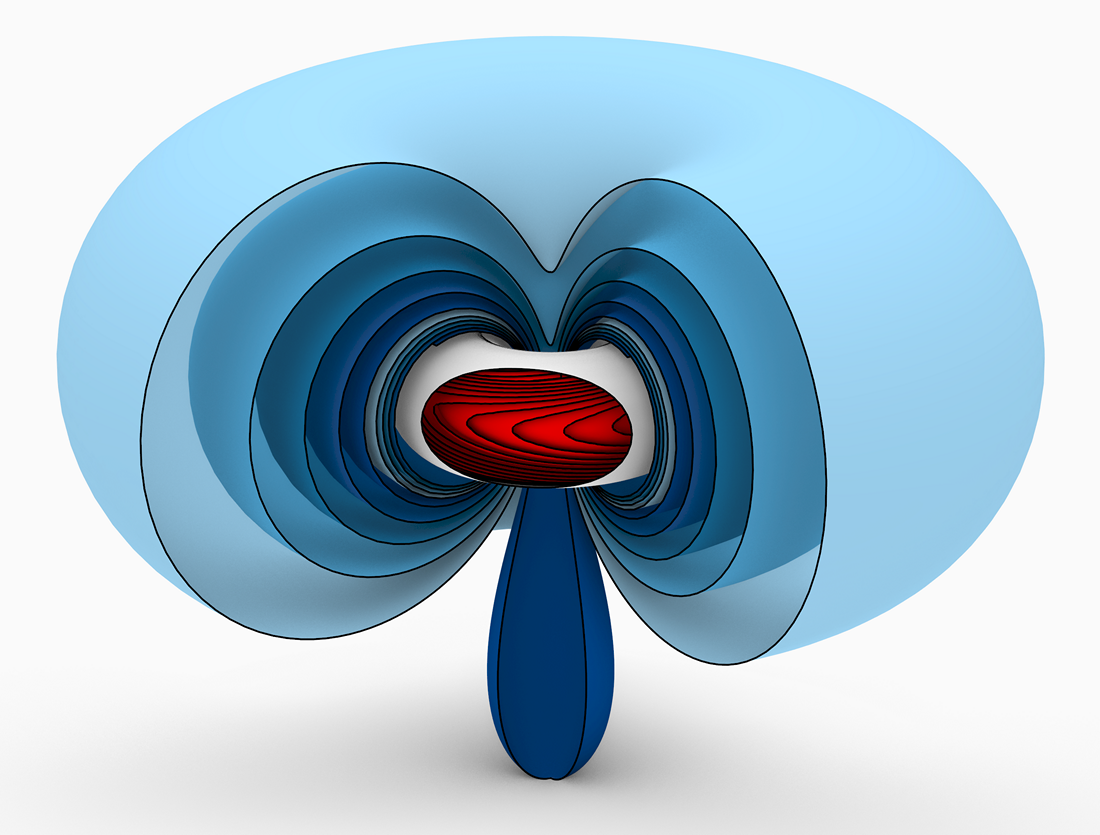}
\end{figure}

The Reeb foliation can be realized in the standard unit sphere of $\mathbb{R}^4$ in such a way that the compact leaf is the Clifford torus and there is an isometry interchanging the two Reeb components, which means that they ``look" the same.
Visualizing the global structure of the Reeb foliation is somehow difficult because it lives on the $3$-sphere which sits in $\mathbb{R}^4$. A quite natural way to do that is projecting it in $\mathbb{R}^3$ via the stereographic projection, missing only one point of the sphere. But the stereographic projection breaks the symmetry of the two Reeb components and therefore one of the Reeb components look easy to grasp and a lot of images of it can be found on textbooks and on the Web, whereas the other component looks unusual. Nevertheless, using the stereographic projection it is possible to visualize the global behaviour of the  foliation. A rendered image of the foliation projected in $\mathbb{R}^3$ is shown in Figure \ref{first rendering 1}.

Such rendering obviously gives an insight of the behaviour of the foliation but very often, especially on a seminar or when dealing with students during a lecture, it could be very useful to have a $3$-dimensional physical model representing this mathematical object.

In general, creating 3D models of mathematical objects has the aim of enhancing understanding, visualization, and engagement with complex concepts, particularly in educational settings. By physically manipulating and exploring these models, students can strengthen their spatial reasoning, observe geometric and topological features more intuitively and gain a deeper grasp of the relationships between different mathematical objects.

Furthermore, the creative process of creating 3D models, coupled with the tangible results, can significantly increase student motivation and engagement in mathematics. Also, 3D models can be used in various educational contexts, from visualizing geometric shapes to understanding complex algorithms or mathematical models of real-world phenomena.

In essence, the motivation for creating 3D models of mathematical objects is to make mathematics more accessible, engaging, and relevant for learners, fostering a deeper understanding and appreciation of its power and applicability.

Coming back to the Reeb foliation, a 3D model of it certainly gives a deeper insight than a rendering. Obviously, the construction of the model has to deal with several problems because of its ``phisycality" but the resolution of these problems is in some sense part of the game.

Our model uses 3D-printing techniques applied to our mathematical context, which will be detailed in Section \ref{rendering}. To our knowledge, the model we have realized is the first attempt of making a physical model of the whole Reeb foliation. It was shown for the first time in September 2024 at the {\it Differential Geometry Workshop} in Brest (France).

In Section \ref{prel}, we will give the basic notions about the theory of foliations and in Section \ref{reeb} we will describe the Reeb foliation of the sphere. Sections \ref{exist models} and \ref{model} are devoted to the discussion of existing models (in fact there are none!) and to our model. In the last section we briefly discuss some future works.
%
%
%
%
%
%
%
%
%
%
\section{Preliminaries}\label{prel}

Let $M$ be a smooth manifold of dimension $n$. Roughly speaking, a \emph{foliation} of dimension $k$ (or {\it codimension} $n-k$) of $M $ is a decomposition of $M$ into a disjoint union of connected submanifolds, called \emph{leaves}, satisfying the following conditions:
\begin{enumerate}
    \item Each leaf $L \subset M $ is a connected immersed submanifold of $M$ of dimension $k$, where $0 \leq k \leq n$.
    \item The foliation is locally trivial: for every point $p \in M$, there exists a coordinate chart $(U, \varphi)$, where $U$ is an open subset of $M$ containing $p$, and $\varphi: U \to \mathbb{R}^n $, such that the leaves in $U$ are given by the preimages
    $$
    \varphi^{-1}(\mathbb{R}^k \times \{c\}),
    $$
    where $c \in \mathbb{R}^{n-k} $ is a constant.
\end{enumerate}
The second condition says that, at least locally, a foliation looks like a $k$-dimensional {\it lasagna}.

\begin{figure}
\caption{Transition maps.}
\centering
\begin{overpic}[width=0.5\textwidth]{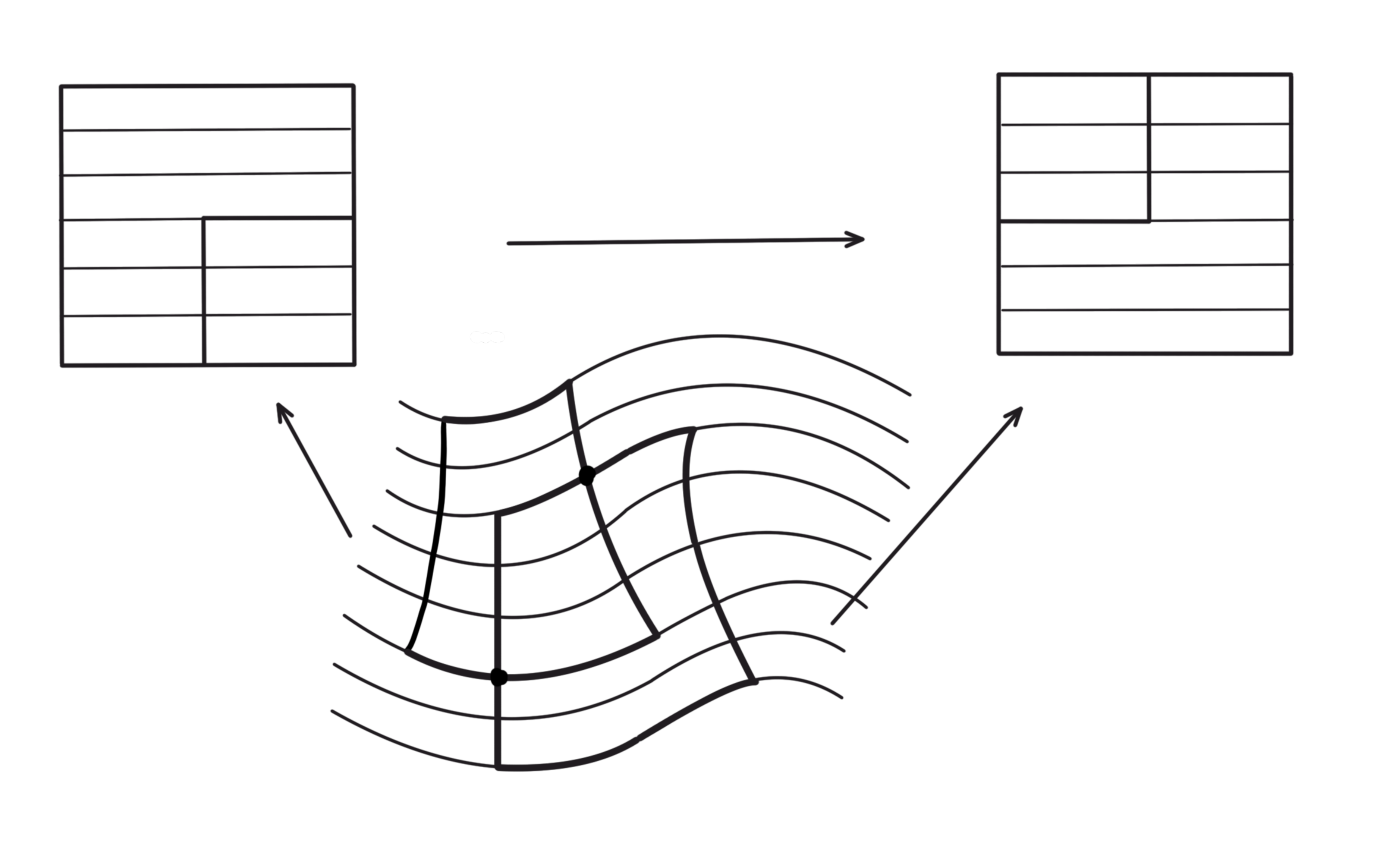}\label{chart diagram}
 \put (40,50) {$\varphi_\beta \circ \varphi_\alpha^{-1} $}
 \put (15,25) {$\varphi_\alpha$}
 \put (73,25) {$\varphi_\beta$}
\end{overpic}
\end{figure}

In a more precise way, a foliation can be described by a  \emph{foliated atlas} on \( M \), which provides the local structure of a foliation by specifying how the leaves of the foliation are represented in coordinate charts. It is a refinement of the concept of a smooth atlas, designed to respect the decomposition of $M $ into leaves.
\begin{Definition}\cite{CC}
A foliated atlas $\mathcal{A} = \{(U_\alpha, \varphi_\alpha)\} $ consists of a collection of charts that cover $M $, where each chart satisfies:
\begin{enumerate}
    \item $U_\alpha \subset M $ is an open subset, and $\varphi_\alpha: U_\alpha \to \mathbb{R}^n $ is a diffeomorphism onto its image.

    \item If $U_\alpha \cap U_\beta \neq \emptyset $, the transition maps $\varphi_\beta \circ \varphi_\alpha^{-1} $ preserve the product structure of the foliation, meaning they are of the form (see Figure \ref{chart diagram}):
$$
\varphi_\beta \circ \varphi_\alpha^{-1}(x, y) = (g(x, y), h(y)),
$$
where $x \in \mathbb{R}^k$ and $y \in \mathbb{R}^{n-k}$.
\end{enumerate}
\end{Definition}

\begin{figure}[h]
\caption{Distinguished chart with plaques.}\label{plaques}
\centering
\includegraphics[width=0.3\textwidth]{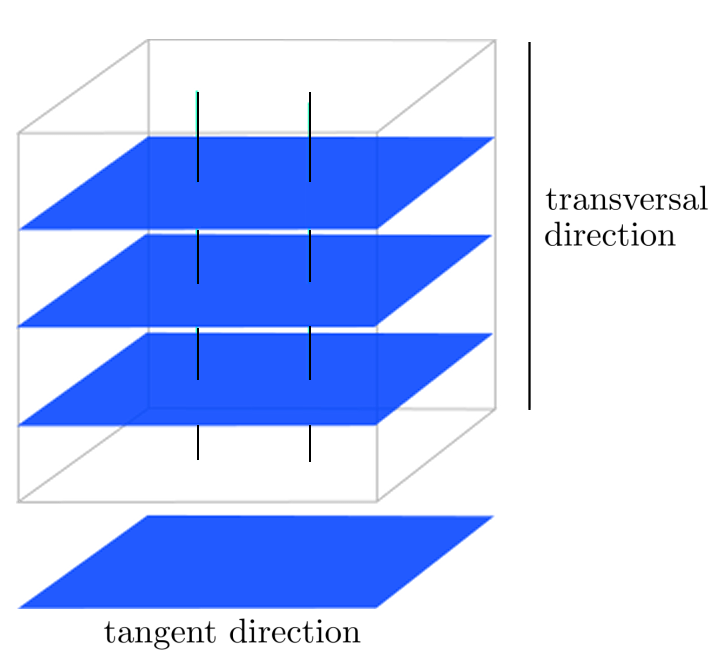}
\end{figure}

To recover the leaves of the foliation from the foliated atlas, one has to consider the {\it plaques} (see Figure~\ref{plaques}) of a chart $\{(U_\alpha, \varphi_\alpha)\}$, which are defined as the connected components of the sets

$$
 \varphi_\alpha^{-1}(\mathbb{R}^k \times \{c\}),
$$
where $\mathbb{R}^k \times \{c\} \subset \mathbb{R}^n$ and $c \in \mathbb{R}^{n-k}$ is constant.

The leaves are therefore given by gluing together the plaques which intersect. More precisely, one defines an equivalence relation on the manifold $M$, by saying that two points $p,q\in M$ are equivalent if there is {\it path of plaques} between them. A path of plaques is a finite number of plaques, say $(P_1, \cdots, P_k)$, such that $p\in P_1$, $q\in P_k$, and each plaque $P_j$ has non-empty intersection with the following one $P_{j+1}$. The resulting equivalence classes are the leaves of the foliation.

Observe that, in general, the leaves of a foliation are immersed submanifolds, not necessarily embedded.

\subsection{Some examples}
The simplest example of foliation (which in fact serves as a local model) is given by the decomposition  of $ \mathbb{R}^n$ into parallel $k$-planes, that is the level sets of the function $\pi: \mathbb{R}^n \to \mathbb{R}^{n-k}$ defined by
$$
\pi (x_1, \cdots, x_n)=(x_{n-k}, \ldots, x_n).
$$
More interesting examples are given by fibre bundles (i.e. tangent and cotangent bundles, principal bundles etc.). In particular, the Hopf fibration of the $3$-sphere is a foliation (see \cite{L} for an elementary explanation of the Hopf fibration).

\begin{figure}[h]
\caption{Level sets of the function \eqref{level sets}, for $n=1$.}\label{submersion-dim 1}
\centering
\includegraphics[width=0.3\textwidth]{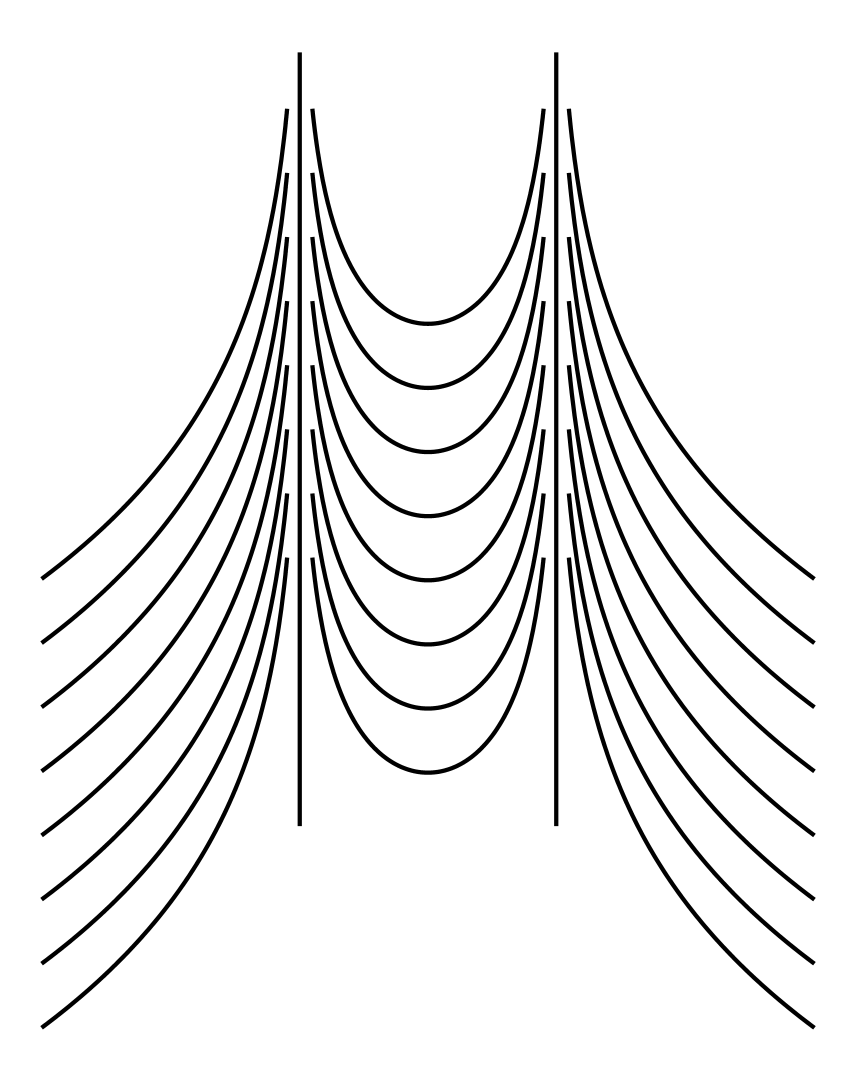}
\end{figure}

\begin{figure}[h]
\caption{Level sets of the function \eqref{level sets}, for $n=2$.}\label{submersion-dim 2}
\centering
\includegraphics[width=0.5\textwidth]{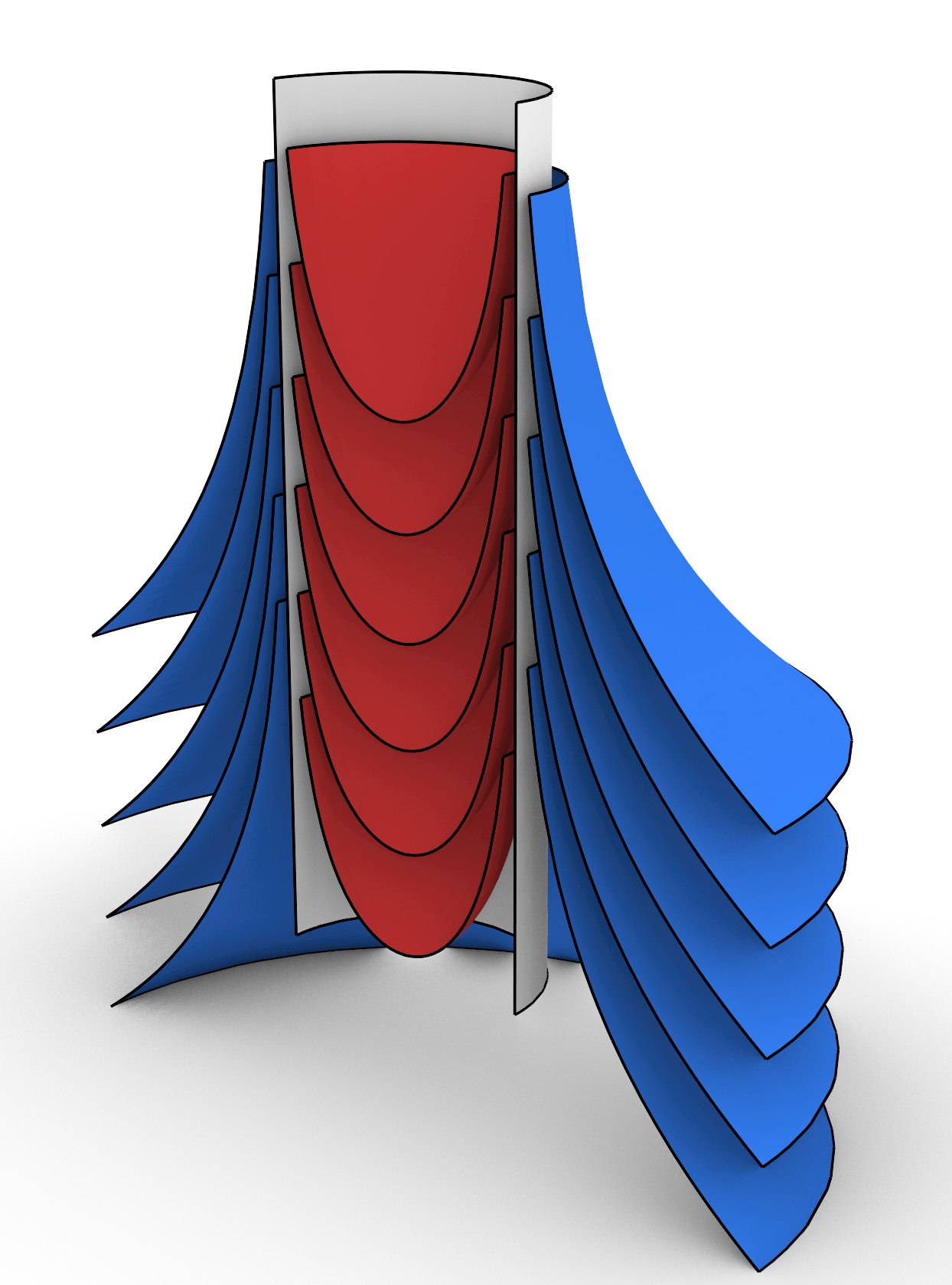}
\end{figure}

More general examples are given by (the connected components of) the level sets of submersions. We give here a simple case which is important in the construction of the Reeb foliation:
\begin{Example}\label{ex-submersion}
Consider the function $f:\mathbb{R}^{n+1} \to \mathbb{R}$ given by
\begin{equation}\label{level sets}
f(x_1, \cdots, x_{n+1})=(r^2-1)\exp{x_{n+1}},
\end{equation}
where $r^2=\sum_1^n  x_i^2$.
A simple calculation shows that $f$ is a submersion and then one obtains a foliation of $\mathbb{R}^{n+1}$ given by the connected components of the level sets of $f$.

For $n=1$ all the leaves are diffeomorphic to $\mathbb{R}$ (Figure~\ref{submersion-dim 1}), but for $n\geq 2$ the leaves are diffeomorphic to  $\mathbb{R}^{n}$ for $r<1$ and to the (generalized) cylinder $\mathbb{S}^{n-1}\times \mathbb{R}$ for $r\geq 2$. In Figure~\ref{submersion-dim 2} are sketched the level sets for $n=2$.
\end{Example}


%


\subsection{The tangent bundle of a foliation and the Frobenius Theorem}
Since the leaves of a $k$-dimensional foliation on $M$ are immersed submanifolds, it is possible to define their tangent spaces and then one can form a subbundle of the tangent bundle $TM$ of $M$, giving rise to a {\it distribution} $\mathcal{D}$. The distribution $\mathcal{D} \subset TM$ satisfies the following condition:
$$
[X, Y] \in \mathcal{D}, \quad \forall X, Y \in \Gamma(\mathcal{D}),
$$
where $\Gamma(\mathcal{D})$ denotes the set of smooth sections of $\mathcal{D}$.

Such a distribution is said to be {\it involutive}. A distribution is said to be {\it completely integrable} if at every point $p$ of $M$ there exists a submanifold passing through $p$ and whose tangent space coincides with the subspace $\mathcal{D}_p$.

The Frobenius theorem provides necessary and sufficient conditions for the integrability of a distribution, which corresponds in fact to the existence of a foliation. A proof of the theorem can be found in \cite{CC}

\begin{Theorem}[Frobenius]
Let $M$ be a smooth manifold, and let  $\mathcal{D}$ be a $k$-dimensional smooth distribution on $M$, that is   $ \mathcal{D} \subset TM $ assigns to each point $p \in M $ a $ k$-dimensional subspace $\mathcal{D}_p \subset T_pM $.

The distribution $\mathcal{D}$ is \emph{completely integrable} (integrable for short) if and only if one of the following equivalent conditions is satisfied:
\begin{enumerate}
\item if $ \mathcal{D} $ is locally spanned by $k$ smooth vector fields $X_1, \dots, X_k$, then $[X_i, X_j] \in \Gamma(\mathcal{D}), \forall i, j=1, \cdots, n$.

\item if $\mathcal{D}$ is locally defined as the kernel of $n-k$ independent $1$-forms $\omega^1, \dots, \omega^{n-k}$, then the following condition holds:
$$
d\omega^i \wedge \omega^1 \wedge \dots \wedge \omega^{n-k} = 0, \quad \forall i=1,\cdots,n-k.
$$

 \item There exists a $k$-dimensional foliation $ \mathcal{F} $ of $M $, such that each leaf $ L$ of $\mathcal{F}$ is an immersed $k$-dimensional submanifold with $T_pL = \mathcal{D}_p$ at every point $p \in L$.

\end{enumerate}

\end{Theorem}

\subsection{Applications of the Frobenius Theorem}

We give here few examples of the application of the Frobenius Theorem.

\begin{enumerate}

\item On $M$, consider the distribution spanned by a non-singular vector field $X$. By the properties of the Lie bracket we have
$$
[X,X]=0
$$
thus the distribution is integrable and the leaves are the integral curves of $X$.

\item On \( \mathbb{R}^3 \), consider the distribution spanned by \( X = \partial_x + y \partial_z \) and \( Y = \partial_y \). Computing the Lie brackets we obtain
$$
[X, Y] = [\partial_x + y \partial_z, \partial_y] = \partial_z.
$$
Since $\partial_z \notin \text{span}(X, Y)$, the distribution is not integrable.

\item In Example \ref{ex-submersion} the distribution is given by the kernel of $df$, which is a non vanishing closed $1$-form and then the distribution is integrable. Of course the leaves of the foliation are the connected components of the level sets of $f$.

\item On $\mathbb{R}^3 $, let $ \mathcal{D} = \ker(\omega)$ where $\omega = dz - xdy+ ydx$. Computing the exterior derivative of $\omega$ gives
$$
d\omega= d(dz - xdy + ydx) = -dx \wedge dy + dy \wedge dx.
$$
Since $d\omega \wedge \omega \neq 0$, we conclude that \( \mathcal{D} \) is not integrable.
\end{enumerate}

\section{The Reeb Foliation}\label{sect-reeb-fol}\label{reeb}
The Reeb foliation of the $3$-sphere $\mathbb{S}^3$ \cite{R} (see also \cite{CC, G}) is constructed by obtaining $\mathbb{S}^3$ as the result of gluing two copies of the solid torus $D^2 \times \mathbb{S}^1$, where $D^2$ is the unit disc and $\mathbb{S}^1$ the unit circle, that is
$$
D^2=\{(x,y)\in \mathbb{R}^2| x^2+y^2\leq 1\} \, ,
$$
$$
\mathbb{S}^1=\{(x,y)\in \mathbb{R}^2| x^2+y^2= 1\} \, .
$$

The gluing of the two copies is made in such a way that the parallel of one copy is glued to a meridian of the other copy ({see \cite{RO} for an explanation of this topological construction}).

A geometric description can be easily visualized by considering $\mathbb{S}^3$ as the $3$-sphere of radius $\sqrt{2}$ centered at the origin of $\mathbb{R}^4$:

$$
\mathbb{S}^3 =\{(x_1,x_2,x_3,x_4)\in \mathbb{R}^4| (x_1)^2+(x_2)^2+(x_3)^2+(x_4)^2=2\}.
$$

Consider the following subsets:
$$
D_1=\{x\in \mathbb{S}^3|(x_1)^2+(x_2)^2\leq 1\} \, , \, D_2=\{x\in \mathbb{S}^3|(x_3)^2+(x_4)^2\leq 1\}.
$$

Then we have $\mathbb{S}^3=D_1\cup D_2$ and $T:=D_1\cap D_2$ is the Clifford Torus (see \cite[p. 140]{DC}).
\\

Each of $D_1$ and $D_2$ is diffeomorphic to $D^2\times \mathbb{S}^1$.
\\
Observe that $D_1$ and $D_2$ are diffeomorphic simply by restriction of the isometry of $\mathbb{R}^4$
$$
G:\mathbb{R}^4 \to \mathbb{R}^4
$$
defined as $G(x_1,x_2,x_3,x_4)=(x_3,x_4,x_1,x_2)$.

An explicit diffeomorphism between $D_1$ and $D^2\times \mathbb{S}^1$ is given by $F:D_1 \to D^2\times \mathbb{S}^1$ defined as follows
$$
F(x_1,x_2,x_3,x_4)=\left( x_1,x_2,\frac{x_3}{\sqrt{(x_3)^2+(x_4)^2}}, \frac{x_4}{\sqrt{(x_3)^2+(x_4)^2}} \right),
$$
whic inverse is given by
$$
F^{-1} (x,y, \cos \alpha, \sin \alpha)=(x,y, (\sqrt{2-x^2-y^2})\cos \alpha,(\sqrt{2-x^2-y^2}) \sin \alpha).
$$

A diffeomorphism between $D_2$ and $D^2\times \mathbb{S}^1$ is given by $F \circ G$.

\begin{Remark}
The isometry $G$ switches meridians and parallels of $T=\partial D_1=\partial D_2$.
\end{Remark}

\subsection*{The Reeb component}
Coming back to Example \ref {ex-submersion}, for $n=2$, we have that $f: \mathbb{R}^3 \to \mathbb{R}$ is given by
$$
f(x,y,z)=(r^2-1) \exp z ,
$$
where $r^2=x^2+y^2$, and the level sets are represented in Figure \ref{submersion-dim 2}.

\begin{figure}[h]
\caption{Level sets of the function \eqref{level sets}, for $n=2$ with $r^2-1\leq0$. The boundary leaf is depicted in white.}\label{reeb component before quotient}
\centering
\includegraphics[width=0.4\textwidth]{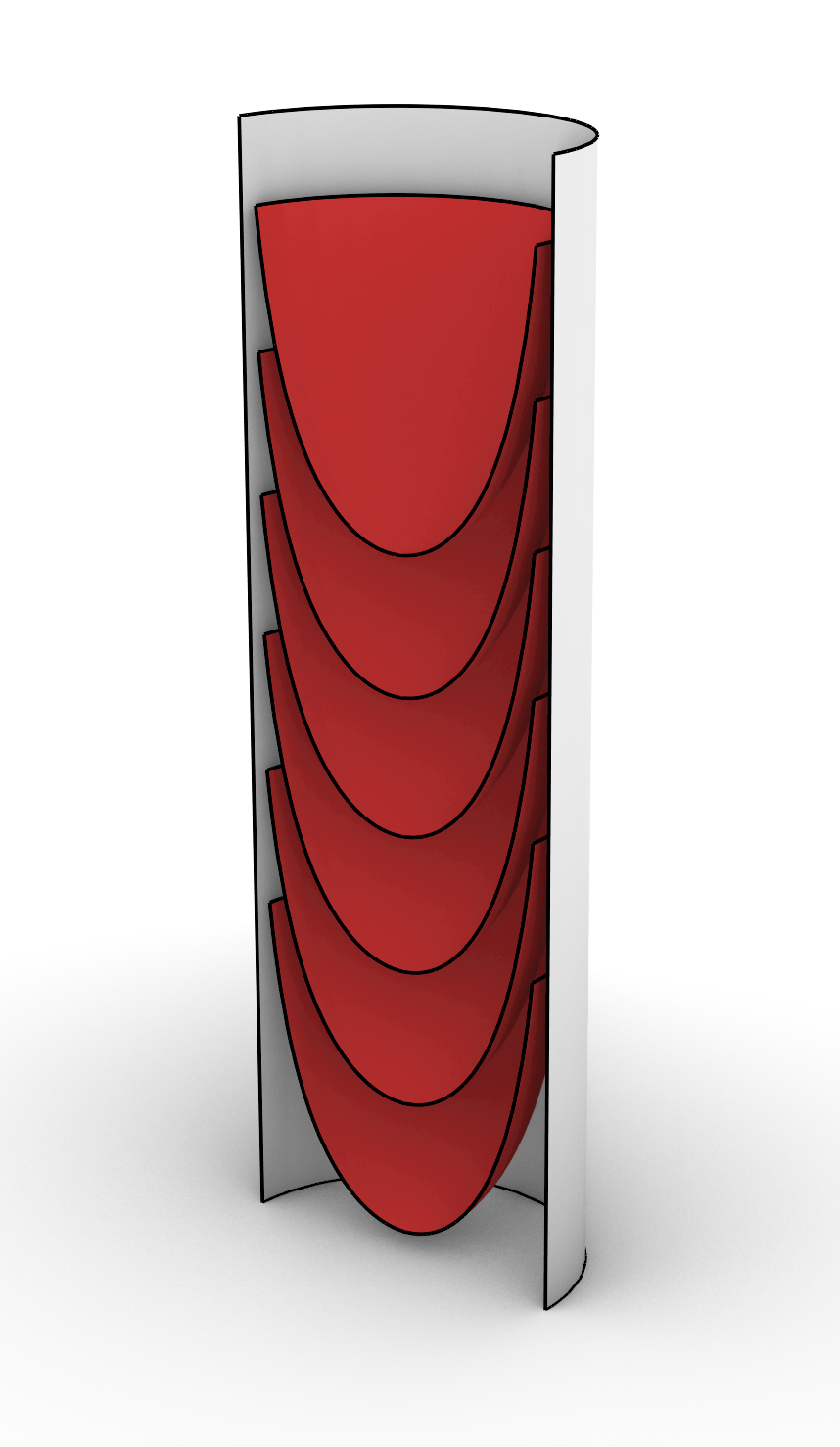}
\end{figure}

If we restrict $f$ to $D^2 \times \mathbb{R}\subset \mathbb{R}^3 $, we obtain a foliation where one leaf is $\mathbb{S}^1\times \mathbb{R}$, the boundary of $D^2 \times \mathbb{R}$ (this is in fact a {\emph {foliation tangent to the boundary} and a precise definition can be found in \cite{CC}). Some leaves of this foliation are represented in Figure~\ref{reeb component before quotient}.

The other leaves are described as follows:
\begin{equation}\label{par-leaves}
\tilde L_c=\left\{ \left( x,y,z=\log \frac{c}{r^2-1} \right) \in \mathbb{R}^3 | r^2-1<0, c\in \mathbb{R}, c<0 \right\}.
\end{equation}
Each leaf $\tilde L_c$ is diffeomorphic to $\mathbb{R}^2$.
\\
Since the foliation is invariant under vertical translations (along the $\mathbb{R}$ factor), one obtains a foliation of the quotient $D^2 \times (\mathbb{R}/\mathbb{Z})=D^2 \times \mathbb{S}^1$, which is called a {\it solid torus}.

To describe explicitly the leaves of $D^2 \times \mathbb{S}^1$ we can just use the complex exponentiation and therefore consider the map
$$
p:D^2 \times \mathbb{R} \to D^2 \times \mathbb{S}^1
$$
given by
$$
p(x,y,z)=(x,y, \exp{2\pi i z}).
$$

The boundary leaf of $D^2 \times \mathbb{R}$ is then projected to the torus $ \mathbb{S}^1\times  \mathbb{S}^1\subset D^2 \times \mathbb{S}^1\subset\mathbb{R}^4$ and the other leaves are still diffeomorphic to $\mathbb{R}^2$, but their behaviour in the solid torus is quite more complicated and can be well visualized in the Figure \ref{reeb component}.

The leaves $L_c:=p(\tilde L_c)$ are given by
\begin{equation}\label{par-leaves-reeb-comp}
\tilde L_c:=\left\{ \left(x,y,\cos {2\pi \log \frac{c}{r^2-1}}, \sin {2\pi \log \frac{c}{r^2-1}} \right)\in \mathbb{R}^4 | r^2=x^2+y^2<1, c\in \mathbb{R}, c<0\right\}
\end{equation}
The manifold $ D^2\times  \mathbb{S}^1$ together with the above defined foliation is called {\it Reeb component}.

\begin{figure}[h]
\caption{Reeb component.}\label{reeb component}
\centering
\includegraphics[width=0.5\textwidth]{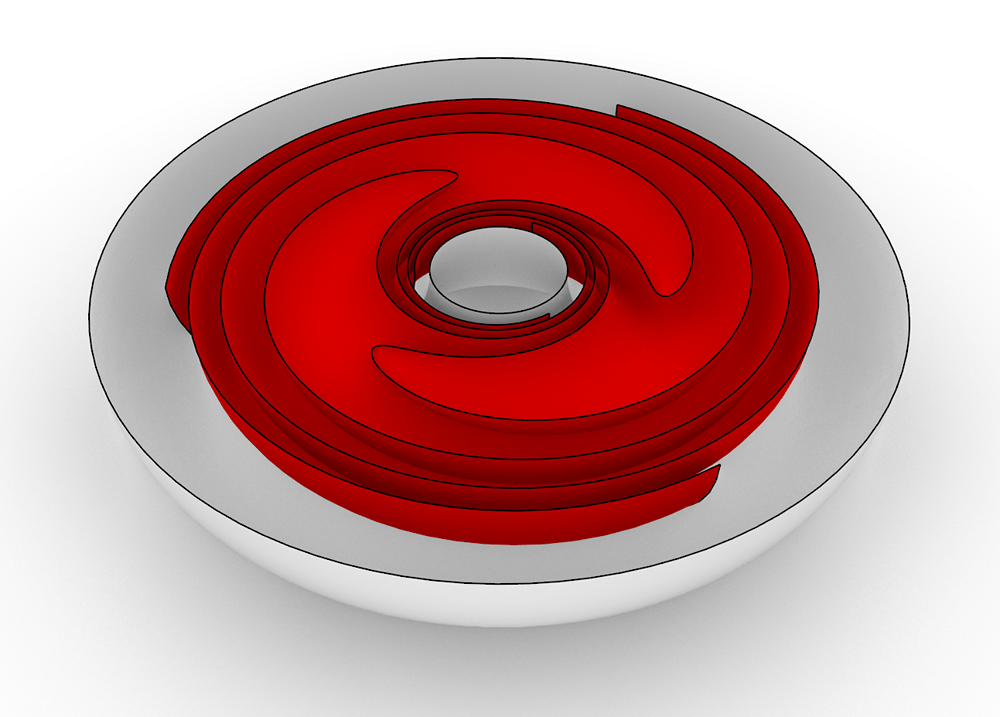}
\end{figure}

\subsection*{The Reeb foliation}
Now recall that the $3$-sphere is obtained by gluing $D_1$ and $D_2$ along their common boundary (the Clifford Torus) as described at the beginning of Section \ref{sect-reeb-fol}.

Since each of $D_1$ and $D_2$ is endowed with a Reeb component via $F^{-1}$ and $(F\circ G)^{-1}$, the Reeb foliation of $\mathbb{S}^3$ is obtained by gluing the two Reeb components along their common boundary. The resulting foliation has one compact leaf (the Clifford Torus) and the other leaves diffeomorphic to $\mathbb{R}^2$.

The non-compact leaves $F^{-1}(\tilde L_c)$ of $D_1$ can be described as follows

\begin{equation}\label{par-leaves-D_1}
\{(x,y,\sqrt{2-r^2}\cos {2\pi \log \frac{c}{r^2-1}},\sqrt{2-r^2} \sin {2\pi \log \frac{c}{r^2-1}})\in \mathbb{R}^4 | r^2=x^2+y^2<1, c\in \mathbb{R}, c<0\}
\end{equation}
and the non-compact leaves $(F\circ G)^{-1}(\tilde L_c)$ of $D_2$ as follows

\begin{equation}\label{par-leaves-D_2}
\left\{ \left(\sqrt{2-r^2}\cos {2\pi \log \frac{c}{r^2-1}},\sqrt{2-r^2} \sin {2\pi \log \frac{c}{r^2-1}},x,y \right)\in \mathbb{R}^4 | r^2=x^2+y^2<1, c\in \mathbb{R}, c<0\right\} .
\end{equation}

\begin{figure}[h]
\caption{Two Reeb components glued along the boundary torus.}\label{reeb component2}
\centering
\includegraphics[width=0.8\textwidth]{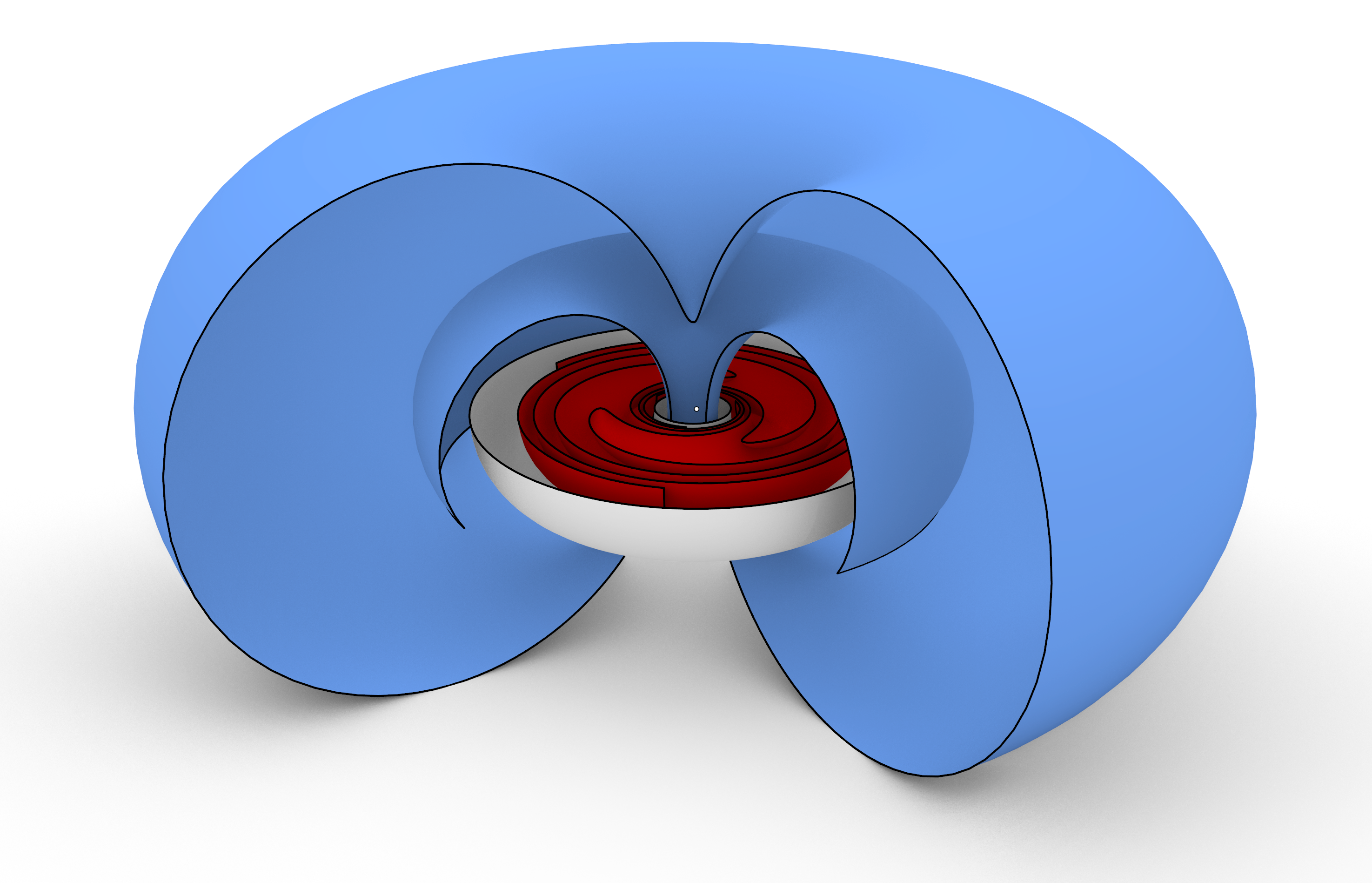}
\end{figure}

\subsection*{Representation in $\mathbb{R}^3$}
In order to visualize the Reeb foliation we need to express it as a set of surfaces lying in ordinary space $\mathbb{R}^3$. This can be done by projecting stereographically our objects from the north pole $N=(0,0,0,\sqrt{2})$, which leads to

 $$
 \begin{aligned}
 \pi_N \colon \mathbb{S}^3\setminus{N} &\to \mathbb{R}^3 \\
 (x_1,x_2,x_3,x_4)   &\mapsto  \left(\frac{\sqrt{2} x_1}{\sqrt{2}-x_4},\frac{\sqrt{2} x_2}{\sqrt{2}-x_4},\frac{\sqrt{2} x_3}{\sqrt{2}-x_4}\right) .
 \end{aligned}
 $$
Via $\pi_N$, the non-compact leaves $F^{-1}(\tilde L_c)$ of $D_1$ are given by

\begin{equation}\label{par-leaves-ster-D_1}
\left\{ \frac{\sqrt{2}}{\sqrt{2}-\sqrt{2-r^2} \sin {2\pi \log \frac{c}{r^2-1}}} (x,y,\sqrt{2-r^2}\cos {2\pi \log \frac{c}{r^2-1}})\in \mathbb{R}^3 | r^2=x^2+y^2<1, c\in \mathbb{R}, c<0\right\}
\end{equation}

and the non-compact leaves $(F\circ G)^{-1}(\tilde L_c)$ of $D_2$ by

\begin{equation}\label{par-leaves-ster-D_2}
\left\{\frac{\sqrt{2}}{\sqrt{2}-y}(\sqrt{2-r^2}\cos {2\pi \log \frac{c}{r^2-1}},\sqrt{2-r^2} \sin {2\pi \log \frac{c}{r^2-1}},x)\in \mathbb{R}^3 | r^2=x^2+y^2<1, c\in \mathbb{R}, c<0\right\} .
\end{equation}

Observe that on $\mathbb{R}^3$ one of the non-compact leaves is (diffeomorphic to) $\mathbb{R}^2 \setminus \pi(N)$.

\section{Rendering and physical models}\label{rendering}

\subsection{Visualization in $\mathbb{R}^3$ via Mathematica\textsuperscript{\textregistered}}
Wolfram  Mathematica\textsuperscript{\textregistered}  is a software system designed for technical computing. It combines symbolic, numerical and graphical capabilities and it is widely used in scientific environment.
In order to plot the leaves of the Reeb foliation in the 3-space we parametrize the non-compact leaves by using polar coordinates in $D^2$. Therefore we obtain the following parametrizations for $\pi_N(F^{-1}(\tilde L_c))$

\begin{equation}\label{par-leaves-ster-polar-D_1}
\phi_c (r,\theta)=\frac{\sqrt{2}}{\sqrt{2}-\sqrt{2-r^2} \sin {2\pi \log \left(\frac{c}{r^2-1}\right)}} \left(r \cos \theta,r \sin \theta,\sqrt{2-r^2}\cos {2\pi \log \left(\frac{c}{r^2-1}\right)}\right)
\end{equation}

We used the \emph{ParametricPlot3D} command to make Mathematica\textsuperscript{\textregistered} draw our objects. Screen drawings produced by Mathematica are polygonal approximations of the surfaces, called \emph{polygonal meshes} or just \emph{meshes}, whose resolution depends on the value of a list of parameters called \emph{PlotPoints}. Once the screen images are satisfactory, we exported the spatial datasets to a format suitable for the needs of rendering and of physical realization. Mathematica\textsuperscript{\textregistered} can export spatial data in a variety of format like for example DXF, DWG, WRML, STL and many others. Among these we have choosen STL (STereo Litography) format, a standard formats used in 3D printing environment. It is a mere list of triangles, given by enumerating the three vertices of each triangle and the three components of the normal to each triangular facet.

\subsection{Physical realization by 3D printing}
In order to use 3D printing techniques to produce physical models of virtual objects, we need to carry on the following steps:
\begin{enumerate}
\item define a polyhedral geometric representation of the object; a triangular net (or mesh) is a good representation;
\item the triangular mesh mentioned on step 1 has to represent a real, solid object, i.e. a polyhedron that encloses a volume. If it does not satisfies this condition, like the object shown in Figure \ref{tubular_neighborhood}, top left, we need to produce a thickened version of it, as shown again in Figure \ref{tubular_neighborhood} on the right;
\item our virtual solid model has to be watertight and must satisfy some specific conditions: it must not have duplicated faces, missing faces, non-manifold edges, misoriented faces;
\item it has to be  expressed in a format manageable by a 3D printer. We used STL format, as said above.
\end{enumerate}
We carry out step 1 by using Mathematica\textsuperscript{\textregistered}.

About step 2: to give a suitable thickness to a surface like in Figure \ref{tubular_neighborhood}, we can proceed in two ways:\\  -with differential geometry tools i.e. by using the normal unit vector field, defined on every regular surface, to generate the so-called \emph{tubular neighborhood} \cite{MH} to our surface (Figure \ref{tubular_neighborhood}) and making Mathematica\textsuperscript{\textregistered} draw directly the thickened object.

-by a CAD program like for example \emph{Rhinoceros 3D} (in \ref{par rhino} we will give some more details about this software), taking the 2D triangular mesh exported from Mathematica\textsuperscript{\textregistered} and applying the tool \emph{Offset Mesh} with the \emph{Solid} option and with a suitable offset distance.
The first approach is more mathematical of course.
We used both of them.

\begin{figure}[h]
\caption{Tubular neighborhood of a regular surface.}\label{tubular_neighborhood}
\centering
\includegraphics[width=0.8\textwidth]{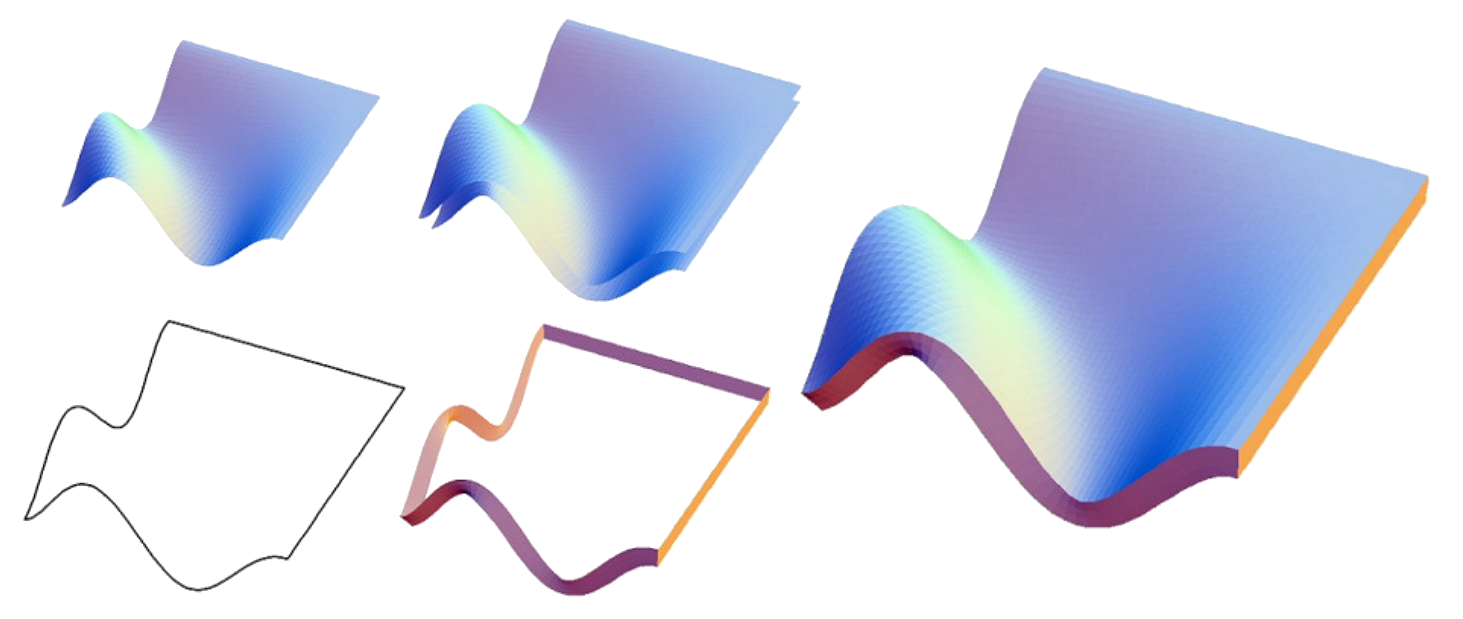}
\end{figure}

To satisfy step 3 means that our representation has to be what in solid modeling is called a \emph{valid solid}. For a comprehensive discussion about the definition of \emph{valid solid} one can refer, for example, to \cite{MA}. We used Rhinoceros to check and, in case of problems, repair the defects the mesh can have.

\subsection{About McNeel Rhinoceros 3D}\label{par rhino}
This is a commercial 3D CAD software developed by Robert McNeel \& Associates, known for its accuracy and flexibility and widely used in automotive, jewelry, shoe and ship  industrial environment for creating, editing, analyzing, rendering, animating NURBS \footnote{NURBS stands for Non-Uniform Rational B-Splines. It is a mathematical model used in computer graphics and 3D modeling to represent curves and surfaces. NURBS are a generalization of both B-splines and Bezier curves, offering a powerful way to create smooth and precise shapes that can be scaled, rotated, and deformed without losing quality.} objects. While Rhino3D is mainly a NURBS modeler, it is equipped with a comprehensive suite of tools for creating and managing mesh objects. Rhino can both create meshes from other object types like NURBS surfaces and solids, and import or export mesh data in various formats.

\subsection{Hardware and software}
The machine used to realize the models is a Stratasys Dimension SST 1200es, which implements the FDM technology (Fused Deposition Modeling)  with soluble supports. It is available at the Laboratory of Mathematics of the University of Cagliari.

Software used: Wolfram Mathematica\textsuperscript{\textregistered} v. 13; Rhinoceros 3D v. 7.

\section{Existing physical models}\label{exist models}
To our knowledge, there are very few images concerning the Reeb foliation. It is possible to find in the literature or on the Web several images of the Reeb component but none of the whole foliation. Clearly, an image like the one in Figure \ref{reeb component2}, where the whole foliation is represented, is more explanatory of the global behaviour than just the one of the Reeb component.

A remarkable image of the Reeb component can be found on \cite{M}.
Such an image seems to suggest the elaboration of a model potentially realizable by 3D printing techniques, but, at the moment, it seems to be only a rendering.


The only attempt to construct of a physical model seems to be the one realized by Bachman (nickname?), found on the Shapeways.com site in march 2025.
At present day (november 2025) that item is not reachable anymore.
Nevertheless, it should be remarked that only a (piece of a) non-compact-leaf of the Reeb component is realized by Bachamn. We can certainly state that our model is the first physical model of the whole Reeb foliation.
\section{The models}\label{model}
We used a filament 3D printer to realize the models described in the present work. The material used is ABS resin (\emph{acrylonitrile butadiene styrene}), a tough and durable plastic resin, and the printer was a Stratasys Dimension SST 1200 thas uses soluble supports.

To explain the models we realized, we will refer to the two families of leaves as {\it inner} and {\it outer} with respect to the image of the Clifford torus under the stereographic projection from $\mathbb{S}^3$ to $\mathbb{R}^3$, which we still call Clifford torus for simplicity, and we use the same terminology for the Reeb components.

The models realized are:
\begin{itemize}

\item A composite model of the Clifford torus, in white color, empty, detachable in two parts, cut along the equator and containing some inner leaves  (in red color) disposable inside the torus as shown in Figure \ref{ogg_realizz_2ab}. The size of the bounding box is $97 \times 97 \times 40 \ mm$). The model represents the inner Reeb component and corresponds to the classical image one finds on texbooks and on the Web.

\item A composite model including an outer leave (the big one in blue color) and a solid torus, detachable in two parts, cut along a vertical plane passing by the symmetry center and joinable magnetically. It can be inserted in the mentioned outer leave as shown in Figure \ref{ogg_realizz_1ab}. The size of the bounding box is $234 \times 234 \times 115 \ mm$. This model represents the outer Reeb component, much less known and  familiar than the previous one.

\item A composite model of a section of the whole foliation, which shows  the Clifford torus (white color) and one sheet for each of the Reeb components: an outer sheet in blue color and a single inner sheet in red color, of which five portions are shown. See Figure \ref{ogg_realizz_3ab}. The dimension of the bounding box is $234 \times 115 \times 83 \ mm$. This model is particularly useful to understand how both the inner and outer leaves accumulate to the boundary Clifford torus.
\end{itemize}

In Figure \ref{allmodels_ab}A one can see the two Reeb components together which form the model to be exposed. Finally, Figure \ref{allmodels_ab}B shows the three models all together. If our aim has been achieved, than looking at the pictures should be more explanatory than any further comment!


\begin{figure}[h]
\caption{Model of the Clifford torus (white) and some of the inner leaves (red).} \label{ogg_realizz_2ab}
\centerline{
\mbox{\includegraphics[width=.4\textwidth]{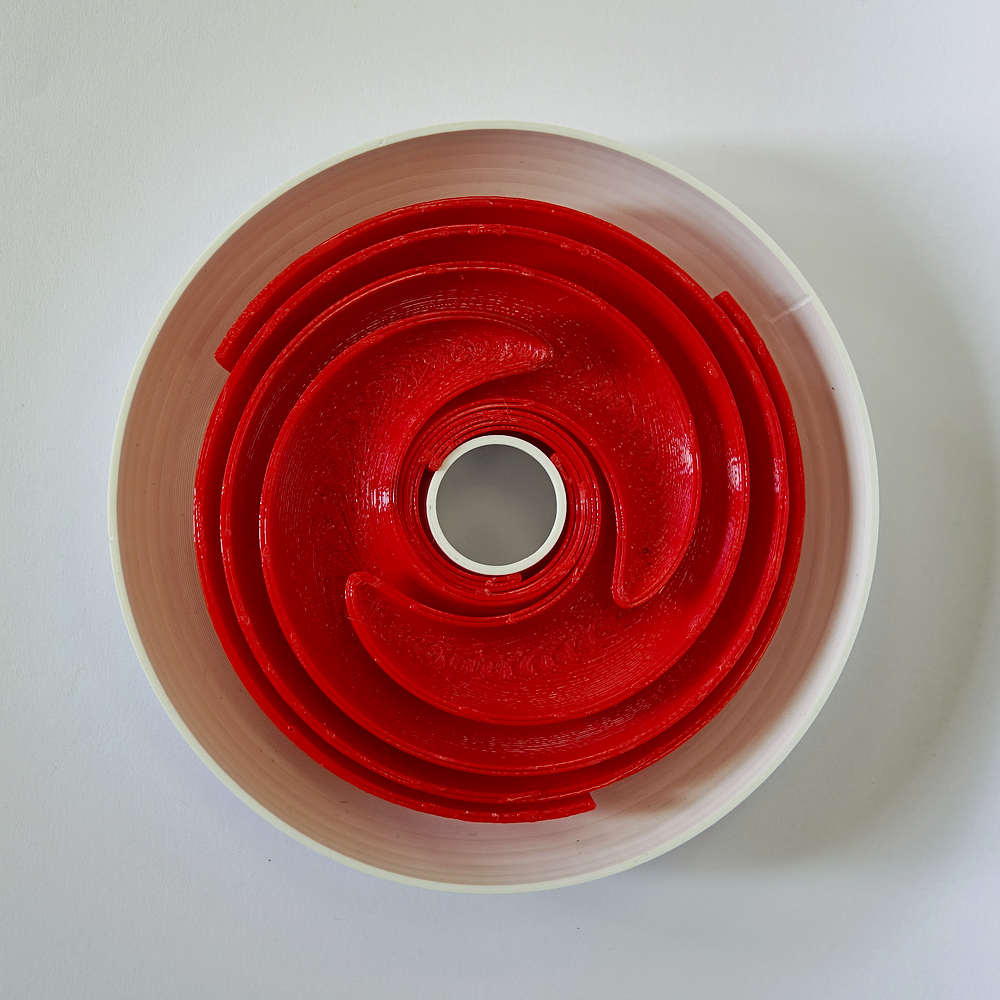}}
\ \
\mbox{\includegraphics[width=.6\textwidth]{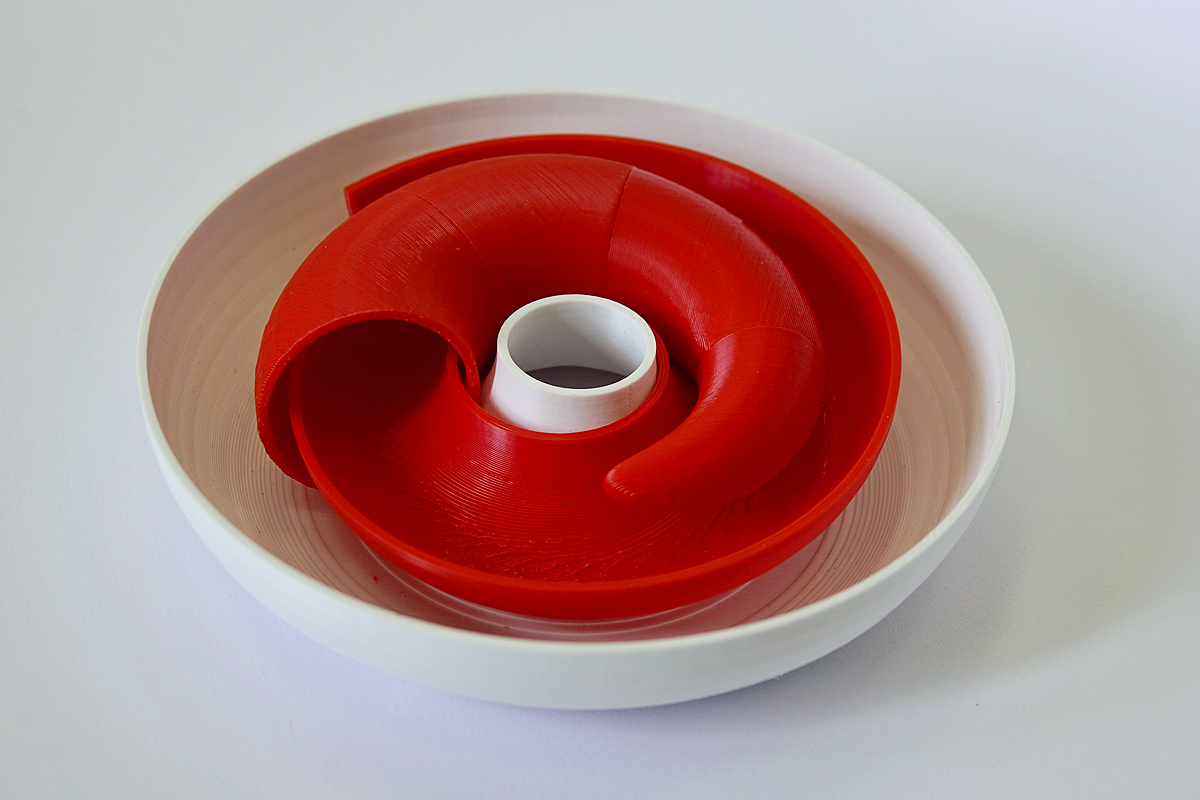}}
}
\end{figure}

\begin{figure}[h]
\caption{Model of the outer leave (blue) and the magnetically detachable solid torus (white).}\label{ogg_realizz_1ab}
\centering
\mbox{\includegraphics[width=.48\textwidth]{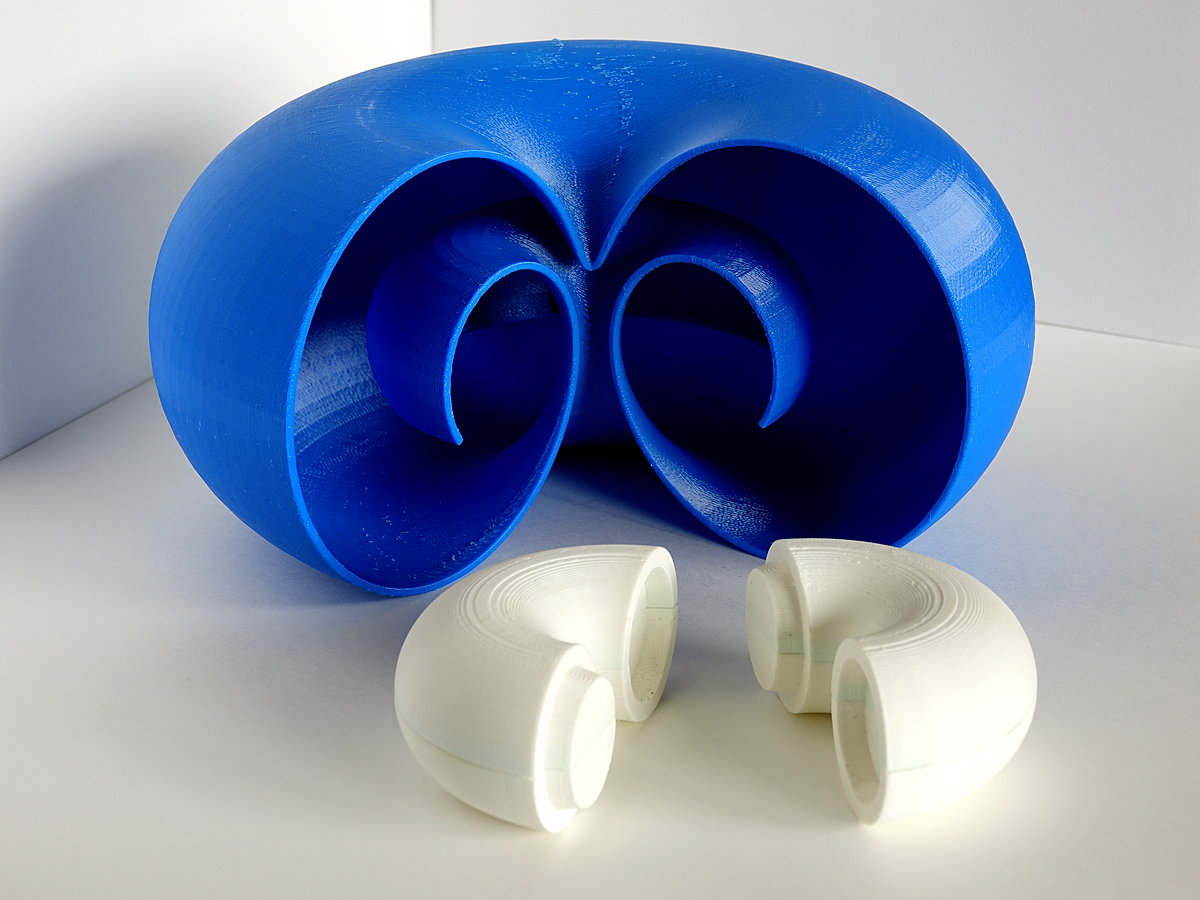}}
\ \
\mbox{\includegraphics[width=.48\textwidth]{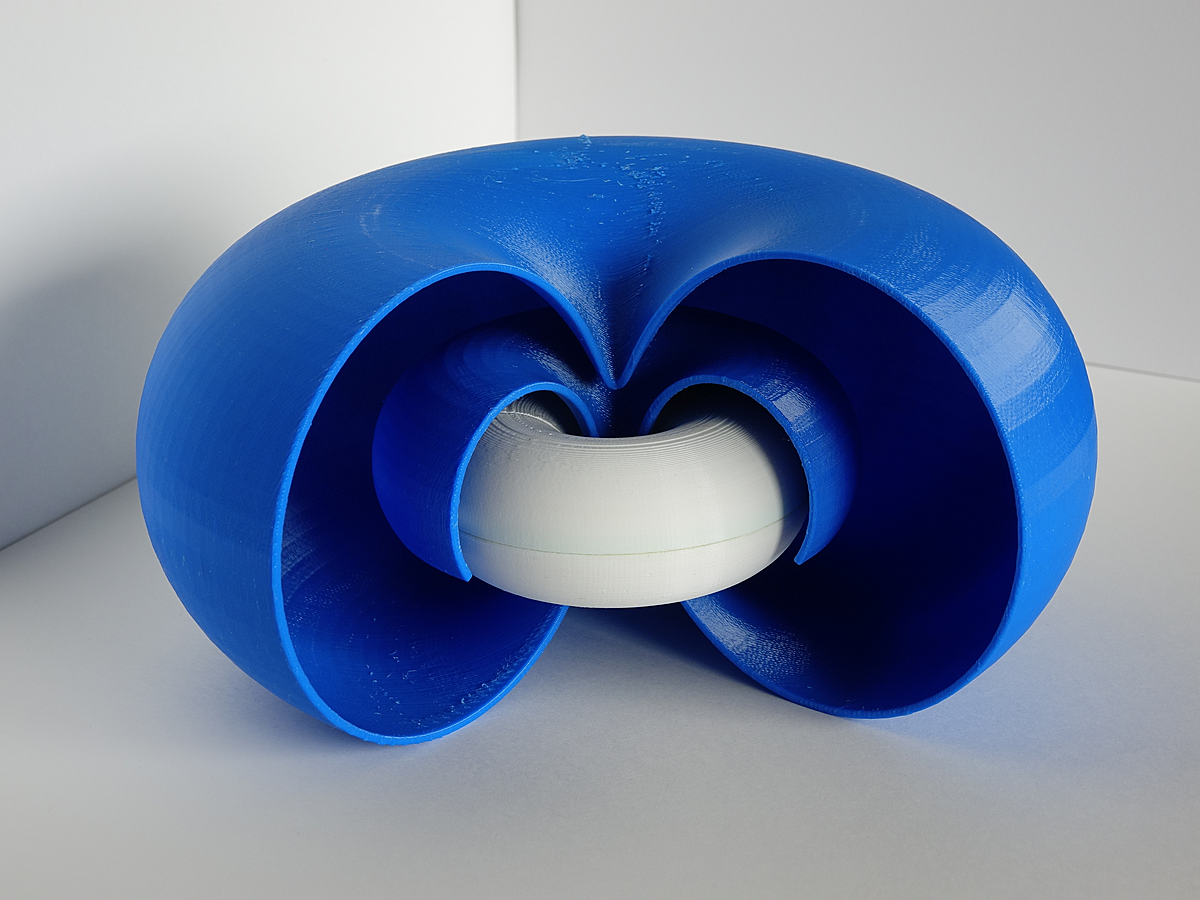}}
\end{figure}

\begin{figure}[h]
\caption{Model of a section of the Reeb foliation showing the Clifford torus and one sheet per Reeb component).} \label{ogg_realizz_3ab}
\centerline{
\mbox{\includegraphics[width=.5\textwidth]{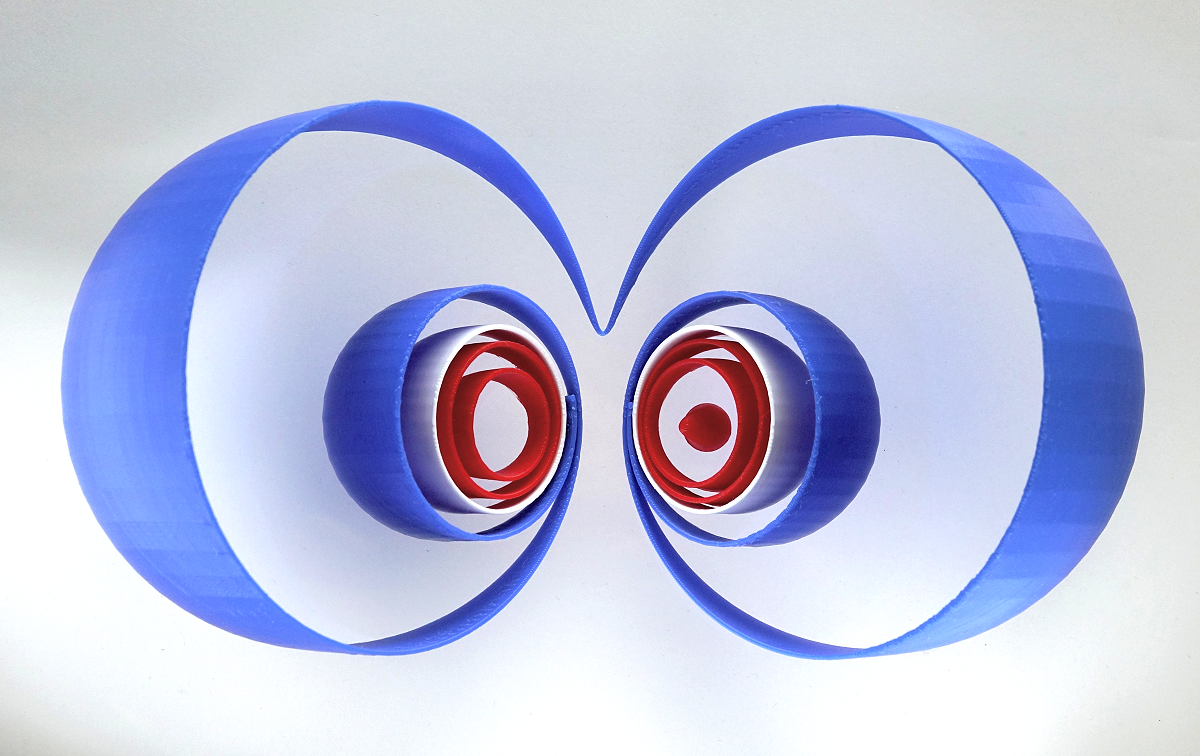}}
\ \
\mbox{\includegraphics[width=.475\textwidth]{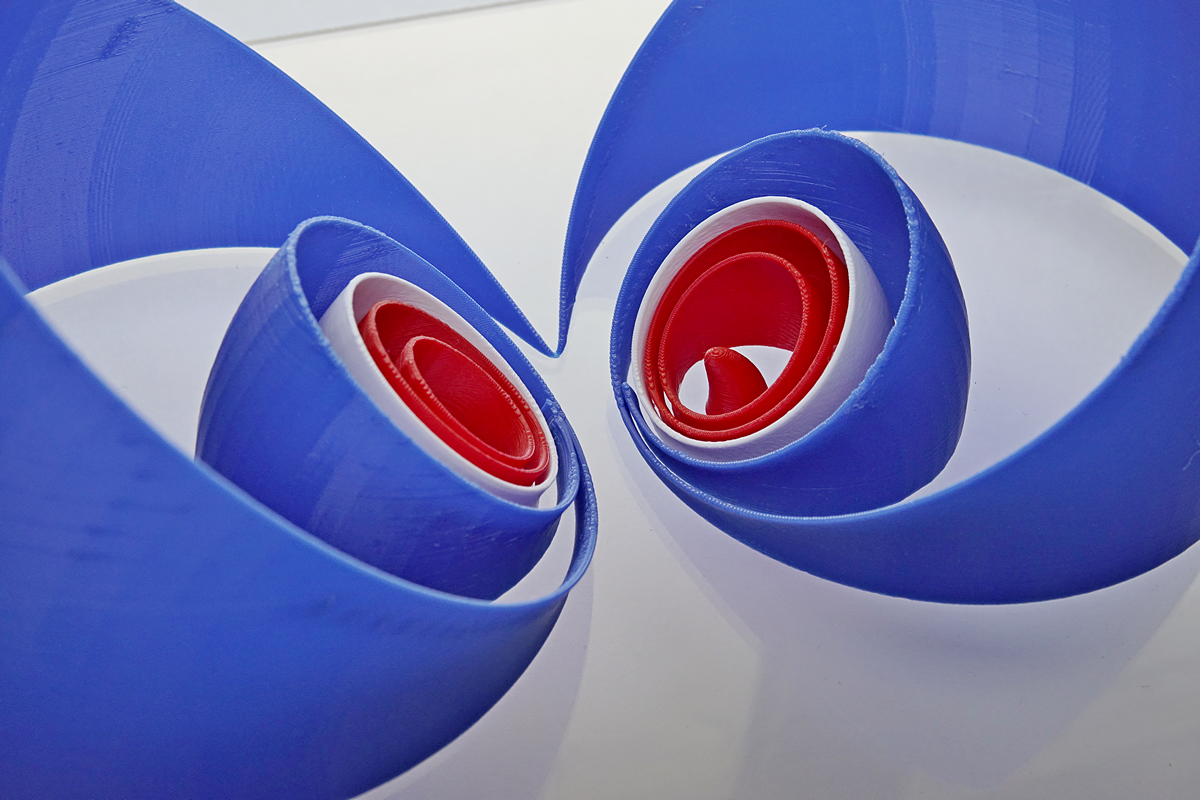}}
}
\end{figure}

\begin{figure}[h]
\caption{On top: the two Reeb components together, which form the main model. Bottom: the same two Reeb components, together with the sectioned foliation.}\label{allmodels_ab}
\vspace{2mm}
\includegraphics[width=0.8\textwidth]{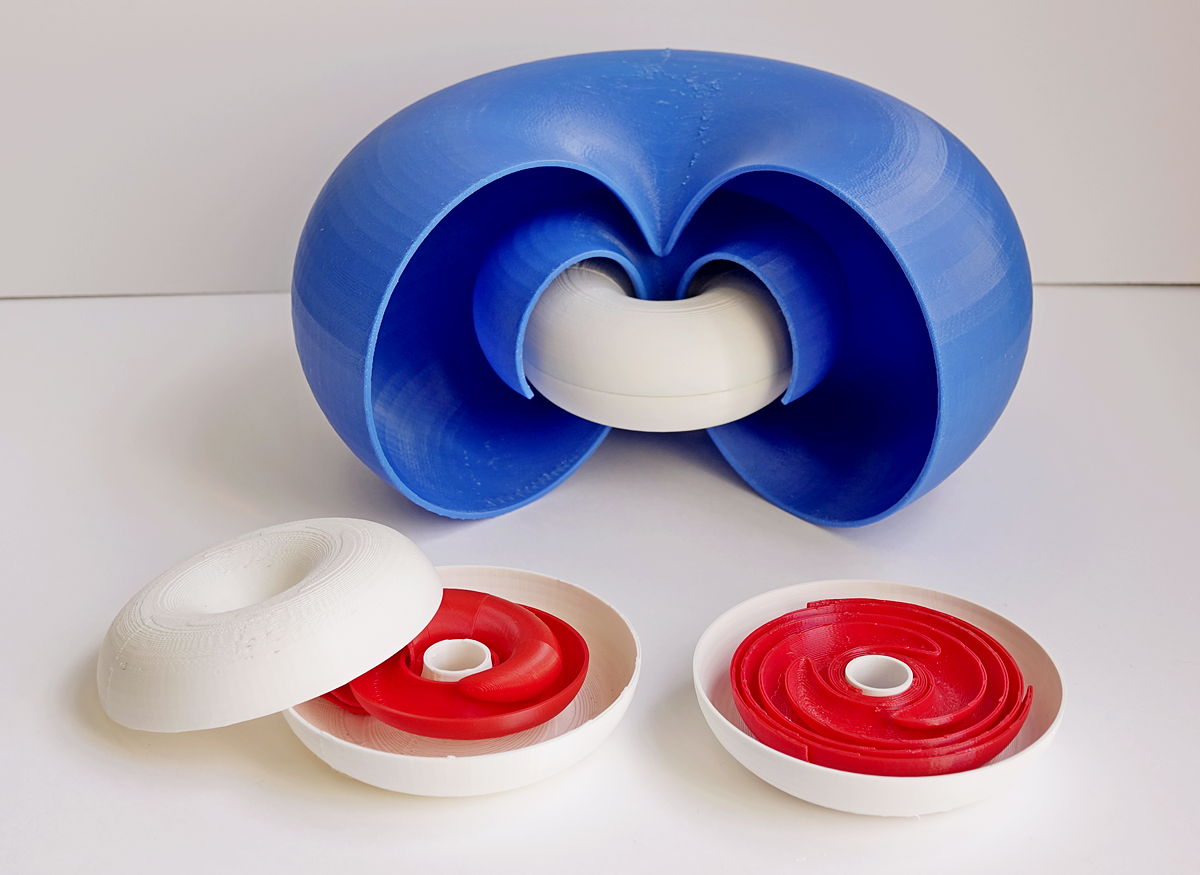}
\vspace{2mm}

\includegraphics[width=0.8\textwidth]{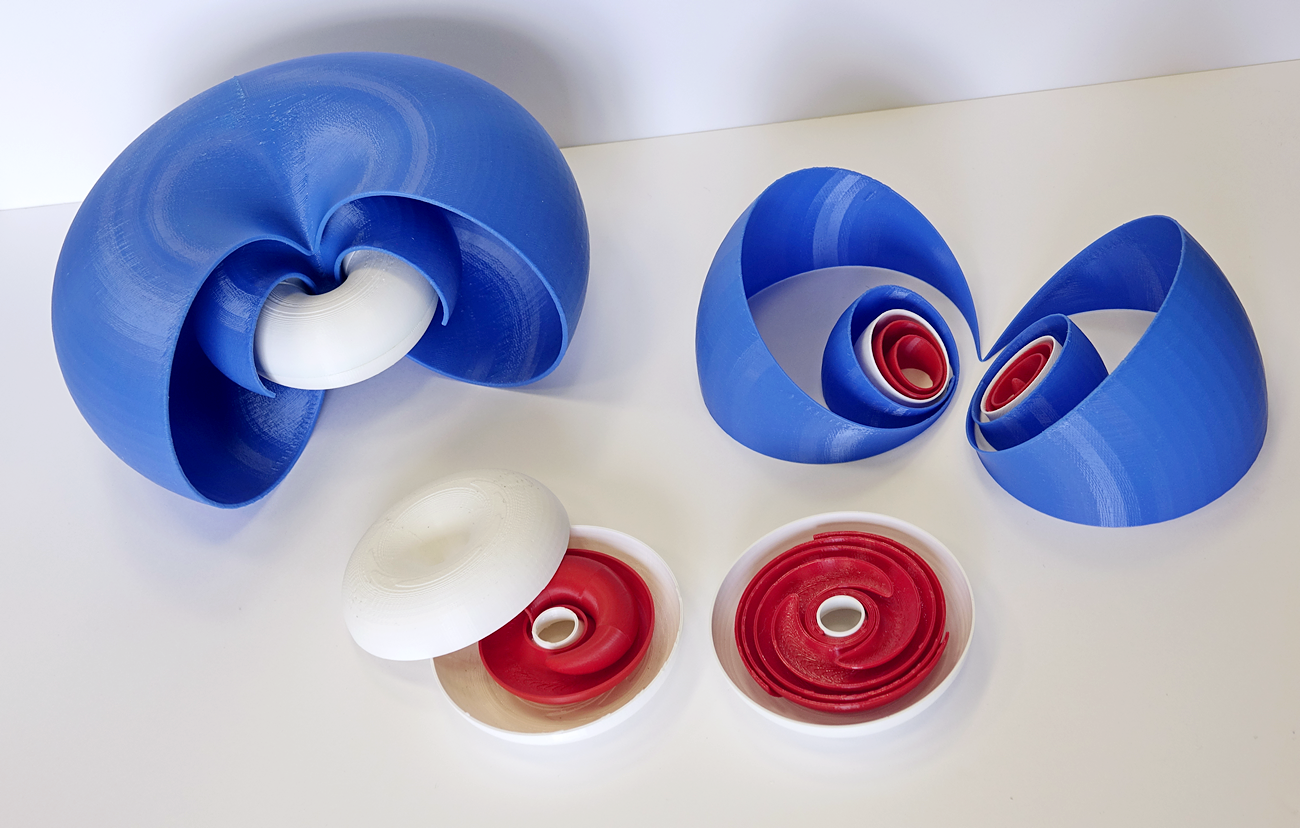}
\end{figure}

\section{Further developments}
To our taste, the model we realized is quite satisfactory, except for the fact that the model of the Reeb component (Figure \ref{ogg_realizz_2ab}) cannot be inserted in the model of the external leaf (the blue one in Figure \ref{ogg_realizz_1ab}
). As an idea for a future work, we would like to design and realize a model completely detachable, where the external leaf would be splittable in two parts to make it possible to place the Reeb component  inside it. The white torus of the Reeb component should be splittable horizontally in two parts joinable magnetically. In this way it should be possible to have a ``self-contained" Reeb foliation reducing the exposition space and enhancing its portability.

\section*{Aknowledgement}
The first author is member of the Italian National Group G.N.S.A.G.A. of INdAM.

This work is partially supported by the project ``GOACT" - CUP: F75F21001210007 - Fondazione di Sardegna - 2020 and by the project ``DAMPAI" - Fondo Crescita Sostenibile - B29J24001610005 - F/350361/01-02/X60.

\end{document}